\documentclass[draft,12pt]{article}

\usepackage{textcomp}

\usepackage[greek,english]{babel} 
\usepackage{amssymb,amsmath} 
\newcommand{\eps}{\varepsilon}
\usepackage{times}

\newtheorem{theorem}{Theorem}

\newtheorem{rem}{Remark}

\newtheorem{pro}{Proposition}

\def\qed{\hfill $\square$}

\begin{document}

\begin{center}
{\bf  Lipschitz perturbations of differentiable implicit functions}\\
Oleg Makarenkov (Voronezh State University, Russia)\\
{\it \small omakarenkov@math.vsu.ru}
\end{center}

{\bf Abstract. } Let $y=f(x)$ be a continuously differentiable
implicit function solving the equation $F(x,y)=0$ with
continuously differentiable $F.$ In this paper we show that if
$F_\eps$ is a Lipschitz function such that the Lipschitz constant
of $F_\eps-F$ goes to 0 as  $\eps\to 0$ then the equation
$F_\eps(x,y)=0$ has a Lipschitz solution $y=f_\eps(x)$ such that
the Lipschitz constant of $f_\eps-f$ goes to 0 as $\eps\to 0$
either. As an application we evaluate the length of time intervals
where the right hand parts of some nonautonomous discontinuous
systems of ODEs are continuously differentiable with respect to
state variables. The latter is done as a preparatory step toward
generalizing the second Bogolyubov's theorem for discontinuous
systems.

\

\noindent {\bf 1. Classical implicit function theorem.} The
classical implicit function theorem can be summarized as follows
(see e.g. \cite{kol}, Ch.~X, \S2, Theorems 1 and 2).

\noindent{\bf Theorem.} {\it Let $X,$ $Y,$ $Z$ be Banach spaces,
$x_0\in X,$ $y_0\in Y$ and $r>0.$ Assume that $F:B_r(x_0)\times
B_r(y_0)\to Z$ satisfies the following conditions
\begin{itemize}
\item[1.] $F(x_0,y_0)=0,$
\item[2.] $F$ is continuous in $B_r(x_0)\times B_r(y_0),$
\item[3.] $F$ is continuously differentiable in $B_r(x_0)\times
B_r(y_0)$ and  $F'$ has in $B_r(x_0)\times B_r(y_0)$ a bounded
inverse.
\end{itemize}
Then there exists $\alpha>0$ and $\beta>0$ such that for any $x\in
B_\alpha(x_0)$ the equation
\begin{equation}\label{eq}
  F(x,y)=0
\end{equation}
has a unique solution $y=f(x)$ in $B_\beta(y_0).$ Moreover, $f$ is
differentiable in $B_\alpha(x_0)$ and
$$f'(x)=-\left[F'_y(x,f(x))\right]^{-1}F'_x(x,f(x))$$
for any $x\in B_\alpha(x_0).$ }

\

\noindent {\bf 2. Main result.} To prove our main result
(theorem~2) we first state the following theorem~1 on the
existence of the implicit function. In the case when the function
$F$ is Lipschitz theorem~1 can be derived from Clark's implicit
function theorem \cite{clark}, but we put the proof since it
appears to be much simpler in our particular situation than that
in \cite{clark}.

\begin{theorem}\label{th1}  Let $X,$ $Y,$ $Z$ be Banach
spaces, $x_0\in X,$ $y_0\in Y$ and $r>0.$ Assume that
$F:B_r(x_0)\times B_r(y_0)\to Z$ satisfies the following
conditions
\begin{itemize}
\item[1.] $F(x_0,y_0)=0,$
\item[2.] $F$ is continuous at $(x_0,y_0),$
\item[3.] $F$ is differentiable at $(x_0,y_0)$ and $F'_y(x_0,y_0)$
has a bounded inverse,
\item[4.] $\|F(x,y)-F(x,y_0)-(F(x_0,y)-F(x_0,y_0))\|\le
L_x\|y-y_0\|,$ for any $x\in B_r(x_0),$ $y\in B_r(y_0),$ where
$L_x\to 0$ as $x\to x_0.$
\end{itemize}
Then there exists $\alpha>0$ and $\beta>0$ such that for any $x\in
B_\alpha(x_0)$ the equation
\begin{equation}\label{eq}
  F(x,y)=0
\end{equation}
has a unique solution $y=f(x)$ in $B_\beta(y_0).$
\end{theorem}

\noindent{\bf Proof.} Let $A_x:B_r(y_0)\to B_r(y_0)$ be defined as
follows $A_x(y)=y-\left[F'_y(x_0,y_0)\right]^{-1}F(x,y).$ Clearly,
the equation
\begin{equation}\label{3}
  A_x(y)=y
\end{equation}
is equivalent to (\ref{eq}).

To prove the existence of solutions to (\ref{3}) we apply the
contracting mappings principle. To this end we show that for any
$\beta>0$ sufficiently small there exists $\alpha>0$ such that for
$x\in B_\alpha(x_0)$ the mapping $A_x$ contracts and  it maps the
ball $B_\beta(y_0)$ into itself. First, using assumptions 4 we
evaluate $A_x(y)-A_x(y_0)$ as follows
\begin{eqnarray*}
  \|A_x(y)-A_x(y_0)\|&=&\left\|y-y_0-\left[F'_y(x_0,y_0)\right]^{-1}(F(x_0,y)-F(x_0,y_0))\right.+\\
  & &
  +\left.\left[F'_y(x_0,y_0)\right]^{-1}\left(F(x_0,y)-F(x_0,y_0)-(F(x,y)-F(x,y_0))\right)\right\|\le\\
  &\le &
  \left\|\left[F'_y(x_0,y_0)\right]^{-1}\left(F(x_0,y)-F(x_0,y_0)-F'_y(x_0,y_0)(y-y_0)\right)\right\|+\\
  & & +L_x\left\|\left[F'_y(x_0,y_0)\right]^{-1}\right\|\cdot\|y-y_0\|.
\end{eqnarray*}
Thus, since $F$ is differentiable at $(x_0,y_0)$ and $L_x\to 0$ as
$x\to x_0$ then for a fixed $\beta>0$ the constant $\lambda>0$ can
be chosen sufficiently small so that
$$
  \|A_x(y)-A_x(y_0)\|\le q\|y-y_0\|,\quad{\rm for \ some\ }q<1{\rm\ and\  any\ }x\in
  B_\alpha(x_0),\ y\in B_\beta(y_0).
$$
Let us now evaluate $\|A_x(y_0)-y_0\|.$ We have
\begin{eqnarray*}
  \|A_x(y_0)-y_0\|&\le
  &\left\|\left[F'_y(x_0,y_0)\right]^{-1}\right\|\cdot \|F(x,y_0)\|=\\
  & &=\left\|\left[F'_y(x_0,y_0)\right]^{-1}\right\|\cdot\|F(x,y_0)-F(x_0,y_0)\|.
\end{eqnarray*}
Therefore, we can diminish $\alpha>0$ in such a way that
$$
  \|A_x(y_0)-y_0\|\le\beta(1-q),\quad{\rm for\ any\ }x\in
  B_\alpha(x_0).
$$
Combining the estimation obtained we arrive to
\begin{eqnarray*}
  \|A_x(y)-y_0\|&\le & \|A_x(y)-A_x(y_0)\|+\|A_x(y_0)-y_0\|\le\\
  & & \le q\|y-y_0\|+\beta(1-q)\le q\beta+\beta(1-q)=\beta.
\end{eqnarray*}
Thus, for any $x\in B_\alpha(x_0)$ the map $A_x$ maps the closed
ball $\overline{B}_\beta(y_0)$ into itself and it contracts in
this ball. Therefore, for any $x\in B_\alpha(x_0)$ the map $A_x$
has a unique fixed point $y=f(x)$ in this ball, that implies
$$
  f(x)=f(x)-\left[F'_y(x_0,y_0)\right]^{-1}F(x,f(x))
$$
or, equivalently, $F(x,f(x))=0.$
 \qed

 \

Next theorem is the main result of the paper. It can be derived
also from the Baitukenov's theorem \cite{bai}. But a proof of
\cite{bai} did not appear in the literature and, thus, we found
reasonable to give a proof independent of the Baitukenov's
theorem.

\begin{theorem} Let $T,V,E,Z$ be Banach spaces and $t_0\in T,$ $v_0\in V$
and $\eps_0\in E.$ Assume that $F:B_r(t_0)\times B_r(v_0)\times
B_r(\eps_0)\to Z$ satisfies the following assumptions
\begin{itemize}
\item[(i)] $F(t_0,v_0,\eps_0)=0,$
\item[(ii)] $F$ is continuous at $(t_0,v_0,\eps_0),$
\item[(iii)] $F'_t(t_0,v_0,\eps_0)$ has a bounded inverse,
\item[(iv)] there exists $L_{\eps,v}\to 0$ as $(\eps,v)\to
(\eps_0,v_0)$ such that
$$
  \|F(t_1,v,\eps)-F(t_2,v,\eps)-F(t_1,v_0,\eps_0)+F(t_2,v_0,\eps_0)\|\le
  L_{\eps,v}\|t_1-t_2\|
$$
for any $t_1,t_2\in B_r(t_0),$ $v\in B_r(v_0),$ $\eps\in
B_r(\eps_0).$
\item[(v)] there exists $K>0$ and $L_\eps\to 0$ as $\eps\to 0$
such that
$$
  \|F(t_1,v_1,\eps)-F(t_1,v_1,\eps)-F(t_2,v_2,\eps_0)+F(t_2,v_1,\eps_0)\|\le
  (L_\eps+K\|t_1-t_2\|)\cdot\|v_1-v_2\|
$$
for any $t_1,t_2\in B_r(t_0),$ $v\in B_r(v_0),$ $\eps\in
B_r(\eps_0).$
\item[(vi)] $(t,v)\to F(t,v,\eps_0)$ is continuously differentiable in $B_r(t_0)\times
B_r(v_0),$
\item[(vii)] $F$ is Lipschitz in $B_r(t_0)\times
B_r(v_0)\times B_r(\eps_0),$
\item[(viii)] The Banach space $T$ possesses the following
property: for any $t\in T$ there exists an element $\{t\}$ of $T$
such that $\{t\}^*t=\|t\|.$ Moreover $\{t\}$ is uniformly bounded
whenever $t$ varies in a bounded set.
\end{itemize}

Then there exists $\alpha>0$ and $\beta>0$ such that for any
$(v,\eps)\in B_\alpha(v_0,\eps_0)$ the equation
\begin{equation}\label{ps2}
  F(t,v,\eps)=0
\end{equation}
has a unique solution $t=\theta(v,\eps)$ in $B_\beta(t_0).$
Moreover, for any $\Delta>0$ there exists $\delta>0$ such that
\begin{equation}\label{INE}
   \|\theta(v_1,\eps)-\theta(v_2,\eps)\|\le\left(\left\|\left[F'_t(t_0,v_0,\eps_0)\right]^{-1}F'_{v}(t_0,v_0,\eps_0)\right\|+\Delta\right)\|v_1-v_2\|
\end{equation}
for any $v_1,v_2\in B_\delta(v_0),$ $\eps\in B_\delta(\eps_0).$
\end{theorem}

\noindent{\bf Proof.} Assumptions (i), (ii), (iii) and (iv) imply
assumptions 1, 2, 3 and 4 of theorem~\ref{th1} with $X=T\times V,$
$x=(v,\eps).$ Therefore, the conclusion about the existence of
$t=\theta(v,\eps)$ solving (\ref{ps2}) follows from
lemma~\ref{th1} and it remains to prove (\ref{INE}).

Let $\Delta>0.$ From the classical implicit function theorem and
assumption (vi) we have that there exists $\delta>0$ such that
$$
  \|\theta(v_1,\eps_0)-\theta(v_2,\eps_0)\|\le\left(-\left[F'_t(t_0,v_0,\eps_0)\right]^{-1}F'_v(t_0,v_0,\eps_0)+\frac{\Delta}{2}\right)\|v_1-v_2\|
$$
for any $v_1,v_2\in B_\delta(v_0).$ Therefore, to prove
(\ref{INE}) it is enough to show that $\delta>0$ can be diminished
in such a way that
$$
  \|\theta(v_1,\eps)-\theta(v_2,\eps)-(\theta(v_1,\eps_0)-\theta(v_2,\eps_0))\|\le\dfrac{\Delta}{2}\|v_1-v_2\|
$$
for any $v_1,v_2\in B_\delta(v_0),$ $\eps\in B_\delta(\eps_0).$

Let $\eta>0$ be fixed. Then by (vi) there exists $d>0$ such that
$$
  F'_t(t_2,v_0,\eps_0)(t_1-t_2)=F(t_1,v_0,\eps_0)-F(t_2,v_0,\eps_0)+\widetilde{\gamma}(t_1,t_2)\cdot\|t_1-t_2\|,
$$
where
\begin{equation}\label{K1}
  \|\widetilde{\gamma}(t_1,t_2)\|\le\eta\quad{\rm for\ any\ }t_1,t_2\in
  B_d(t_0).
\end{equation}
Without loss of generality we can assume that
\begin{equation}\label{ASSU}
  0<d<\eta.
\end{equation}
 For an auxiliary  $v\in B_\delta(v_0)$ we consider
\begin{eqnarray*}
 t_1-t_2&=&\left[F'_t(t_2,v_0,\eps_0)\right]^{-1}F'_t(t_2,v_0,\eps_0)(t_1-t_2)=\nonumber\\
 & =&
 \left[F'_t(t_2,v_0,\eps_0)\right]^{-1}\left(F(t_1,v_0,\eps_0)-F(t_2,v_0,\eps_0)+\widetilde{\gamma}(t_1,t_2)\cdot\|t_1-t_2\|\right)=\nonumber\\
 &=&\left[F'_t(t_2,v_0,\eps_0)\right]^{-1}(F(t_1,v,\eps)-F(t_2,v,\eps)+\widetilde{\gamma}(t_1,t_2)\cdot\|t_1-t_2\|+\nonumber\\
 & & \qquad \qquad \qquad \qquad
 +F(t_1,v_0,\eps_0)-F(t_2,v_0,\eps_0)-(F(t_1,v,\eps)-F(t_2,v,\eps))).
\end{eqnarray*}
By (iv) we can diminish $\delta>0$ in such a way that
$$
  \|F(t_1,v_0,\eps_0)-F(t_2,v_0,\eps_0)-(F(t_1,v,\eps)-F(t_2,v,\eps_0))\|\le\eta\|t_1-t_2\|
$$
for any $t_1,t_2\in B_d(t_0),$ $v\in B_\delta(v_0)$ and $\eps\in
B_\delta(\eps_0).$ Therefore, taking into account (\ref{K1}) we
have that
\begin{equation}\label{PR}
  t_1-t_2=\left[F'_t(t_2,v_0,\eps_0)\right]^{-1}(F(t_1,v,\eps)-F(t_2,v,\eps)-\widehat{\gamma}(t_1,t_2,v)\|t_1-t_2\|),
\end{equation}
where $\|\widehat{\gamma}(t_1,t_2,v)\|\le 2\eta,$ for any
$t_1,t_2\in B_d(t_0),$ $v\in B_\delta(v_0),$ $\eps\in
B_\delta(\eps_0).$

Since
$F(\theta(v_2,\eps),v_2,\eps)=0=F(\theta(v_1,\eps),v_1,\eps)$
taking $t_1=\theta(v_1,\eps),$ $t_2=\theta(v_2,\eps),$ $v=v_2$ we
obtain from (\ref{PR}) that
\begin{eqnarray*}
 \theta(v_1,\eps)-\theta(v_2,\eps)&=&\left[F'_t(\theta(v_2,\eps),v_0,\eps_0)\right]^{-1}(F(\theta(v_1,\eps),v_2,\eps)-
   F(\theta(v_1,\eps),v_1,\eps)+\\
   & & +
   \widehat{\gamma}(\theta(v_1,\eps),\theta(v_2,\eps),v_2)\cdot\|\theta(v_1,\eps)-\theta(v_2,\eps)\|)
\end{eqnarray*}
for any $v_1,v_2\in B_\delta(v_0),$ $\eps\in B_\delta(\eps_0).$ By
assumption (viii) we have
\begin{eqnarray*}
  \theta(v_1,\eps)-\theta(v_2,\eps)&=&\left[F'_t(\theta(v_2,\eps),v_0,\eps_0)\right]^{-1}(F(\theta(v_1,\eps),v_2,\eps)-F(\theta(v_1,\eps),v_1,\eps))+\\
  & &
  +\widehat{\gamma}(\theta(v_1,\eps),\theta(v_2,\eps),v_2)\left\{\overline{\theta(v_1,\eps)-\theta(v_2,\eps)}\right\}
  (\theta(v_1,\eps)-\theta(v_2,\eps)),
\end{eqnarray*}
where $(v_1,v_2,\eps)\mapsto
\left\{\overline{\theta(v_1,\eps)-\theta(v_2,\eps)}\right\}$ is
bounded on $B_\delta(v_0)\times B_\delta(v_0)\times
B_\delta(\eps_0).$ Since $\eta>0$ can be chosen sufficiently small
we can consider that
$I-\widehat{\gamma}(\theta(v_1,\eps),\theta(v_2,\eps),v_2)\left\{\overline{\theta(v_1,\eps)-\theta(v_2,\eps)}\right\}$
is invertible for any $v_1,v_2\in B_\delta(v_0),$ $\eps\in
B_\delta(\eps_0).$ Thus, we can rewrite the previous expression as
follows
\begin{eqnarray*}
  \theta(v_1,\eps)-\theta(v_2,\eps)&=&\left(I-\widehat{\gamma}(\theta(v_1,\eps),\theta(v_2,\eps),v_2)\left\{\overline{\theta(v_1,\eps)-\theta(v_2,\eps)}\right\}\right)^{-1}\circ\\
  & &
  \circ\left[F'_t(\theta(v_2,\eps),v_0,\eps_0)\right]^{-1}(F(\theta(v_1,\eps),v_2,\eps)-F(\theta(v_1,\eps),v_1,\eps))
\end{eqnarray*}
for any $v_1,v_2\in B_\delta(v_0),$ $\eps\in B_\delta(\eps_0).$

By the conclusion of theorem~1 the function $\theta$ is continuous
at $(v_0,\eps_0)$ and we can diminish $\delta>0$ also in such a
way that
\begin{equation}\label{ASSU1}
  \|\theta(v,\eps)-t_0\|\le\frac{d}{2}\quad{\rm for\ any\ }v\in
  B_\delta(v_0),\ \eps\in B_\delta(\eps_0).
\end{equation} Combining (\ref{ASSU}) and (\ref{ASSU1}) we have
\begin{equation}\label{OO}
  \|\theta(v,\eps)-\theta(v,\eps_0)\|\le\eta,\quad{\rm for\ any\
  }v\in B_\delta(v_0),\ \eps\in B_\delta(\eps_0).
\end{equation}
Thus, using assumption (vii) we have the following expression for
$\theta(v_1,\eps)-\theta(v_2,\eps)$
\begin{equation}\label{FO}
  \theta(v_1,\eps)-\theta(v_2,\eps)=\left[F'_t(\theta(v_2,\eps),v_0,\eps_0)\right]^{-1}(F(\theta(v_1,\eps),v_2,\eps)-
  F(\theta(v_1,\eps),v_1,\eps))+\gamma(v_1,v_2,\eps),
\end{equation}
where
$$
  \|\gamma(v_1,v_2,\eps)\|\le L_\eta\|v_1-v_2\|,\quad{\rm for\
  any\ } v_1,v_2\in B_\delta(v_0),\ \eps\in B_\delta(\eps_0){\rm\
  and\ }L_\eta\to 0{\rm\ as\ }\eta\to 0
$$
(of course, $\delta>0$ depends on $\eta>0$ as well).

Formula (\ref{FO}) allows us to evaluate
$\|\theta(v_1,\eps)-\theta(v_2,\eps)-\theta(v_1,\eps_0)+\theta(v_2,\eps_0)\|$
as follows
\begin{eqnarray*}
  & &
  \|\theta(v_1,\eps)-\theta(v_2,\eps)-\theta(v_1,\eps_0)+\theta(v_2,\eps_0)\|\le\\
  & & \qquad
  \le\left\|\left[F'_t(\theta(v_2,\eps_0),v_0,\eps_0)\right]^{-1}\right\|\cdot\\
  & & \qquad\qquad
  \cdot\|F(\theta(v_1,\eps),v_2,\eps)-F(\theta(v_1,\eps),v_1,\eps)-F(\theta(v_1,\eps_0),v_2,\eps_0)+F(\theta(v_1,\eps_0),v_1,\eps_0)\|+\\
  & & \qquad
  +\left\|\left[F'_t(\theta(v_2,\eps),v_0,\eps_0)\right]^{-1}-\left[F'_t(\theta(v_2,\eps_0),v_0,\eps_0)\right]^{-1}\right\|\cdot\\
  & &
  \qquad\qquad\cdot\|F(\theta(v_1,\eps),v_2,\eps)-F(\theta(v_1,\eps),v_1,\eps)\|+\\
  & & \qquad+\|\gamma(v_1,v_2,\eps)\|+\|\gamma(v_1,v_2,\eps_0)\|,
\end{eqnarray*}
for any $v_1,v_2\in B_\delta(v_0),$ $\eps\in B_\delta(\eps_0).$
From (\ref{OO}) we can conclude that
$$
  \left\|\left[F'_t(\theta(v_2,\eps),v_0,\eps_0)\right]^{-1}-\left[F'_t(\theta(v_2,\eps_0),v_0,\eps_0)\right]^{-1}\right\|\le
  K_\eta
$$
for any $v\in B_\delta(v_0),$ $\eps\in B_\delta(\eps_0),$ where
$K_\eta\to 0$ as $\eta\to 0.$ Then, fixing some $K_1>0$ such that
$\left[F'_t(\theta(v_2,\eps_0),v_0,\eps_0)\right]^{-1}\le K_1$ for
any $v_2\in B_\delta(v_0)$ and applying assumptions (v) and (vii)
we have
\begin{eqnarray*}
 & & \|\theta(v_1,\eps)-\theta(v_2,\eps)-\theta(v_1,\eps_0)+\theta(v_2,\eps_0)\|\le\\
 & & \qquad \le
 K_1(L_\eps+K\|\theta(v_1,\eps)-\theta(v_1,\eps_0)\|)\cdot\|v_1-v_2\|+K_\eta
 L\|v_1-v_2\|+2L_\eta\|v_1-v_2\|\le\\
 & & \qquad \le (K_1 L_\eps+K_1 K\eta+K_\eta L+2
 L_\eta)\cdot\|v_1-v_2\|
\end{eqnarray*}
for any $v_1,v_2\in B_\delta(v_0),$ $\eps\in B_\delta(\eps_0).$

Therefore, given $\Delta>0$ we can find $\eta>0$ and $\delta>0$
(which depends on $\eta>0$) such that
$$
  \|\theta(v_1,\eps)-\theta(v_2,\eps)-(\theta(v_1,\eps_0)-\theta(v_2,\eps_0))\|\le\dfrac{\Delta}{2}\|v_1-v_2\|
$$
for any $v_1,v_2\in B_\delta(v_0),$ $\eps\in B_\delta(\eps_0).$
Thus the proof is complete.

\qed

\

\noindent{\bf 3. An application.} Consider the second order
differential equation
\begin{equation}\label{SO}
  \ddot u+u=-\eps{\rm sign}(u)+\eps g(t,u,\dot u),
\end{equation}
where $g$ is continuously differentiable and $2\pi$-periodic in
time. The change of variables
$$
  \left(\begin{array}{c}
  u(t)\\ \dot u(t)\end{array}\right)= \left(\begin{array}{cc}
  \cos t & \sin t\\ -\sin t & \cos t\end{array}\right) \left(\begin{array}{c}
  x_1(t)\\ x_2(t)\end{array}\right)
$$
allows us to transform (\ref{SO}) into the following system
\begin{equation}\label{S1}
  \begin{array}{rcl}
    \dot x_1 &=&\eps\sin(t){\rm sign}(x_1\cos t+x_2\sin t)\\
   \dot x_2 &=&-\eps\cos(t){\rm sign}(x_1\cos t+x_2\sin t)
  \end{array} +\eps\cdot\boxed{\begin{array}{c}{\rm continuously\ differentiable\
  terms}\end{array}}\ .
\end{equation}

We assume that for any $v\in\mathbb{R}^2$ system (\ref{S1}) has an
unique absolutely continuous solution $x(\cdot,v,\eps)$ defined on
$[0,T]$ and such that
\begin{itemize}
\item[(F)] $x(t,v,\eps)$ possesses the following representation
$x(t,v,\eps)=v+\eps y(t,v,\eps),$ where $(t,v)\mapsto y(t,v,\eps)$
is locally Lipschitz uniformly with respect to small $\eps>0.$
\end{itemize}

\begin{rem} A natural example when all the imposed assumptions are satisfied
is when (\ref{SO}) models a pendulum with dry friction (see
\cite[example~2.2.3]{kunze}).
\end{rem}

Consider
$$
  F(t,v,\eps)=x_1(t,v,\eps)\cos t+x_2(t,v,\eps)\sin t.
$$
The following proposition is crucial when generalizing the second
Bogolyubov's theorem \cite{bog} for discontinuous systems of form
(\ref{S1}).
\begin{pro} Assume that (F) is satisfied. Let
$t_0\in(a,b)\subset(0,2\pi)$ be the only zero of $F(\cdot,v_0,0)$
on $[a,b]$ and define $$R=\dfrac{1}{\left|-[v_0]_1\sin
t_0+[v_0]_2\cos t_0\right|}.$$ Then there exists a function
$\theta:B_\delta(v_0)\times B_\delta(v_0)\to [a,b]$ such that
given $\Delta>0$ there exists $\delta>0$ such that $(t,v)\mapsto
F(t,v,\eps)$ does not vanish on
\begin{equation}\label{NV}
 \left([a,b]\backslash
 \left[\theta(v_1,\eps)-(R+\Delta)\|v_1-v_2\|,\theta(v_1,\eps)+(R+\Delta)\|v_1-v_2\|\right]\right)\times[v_1,v_2]
\end{equation}
whenever $v_1,v_2\in B_\delta(v_0),$ $\eps\in (0,\delta).$
\end{pro}

\noindent{\bf Proof.} Let us show that assumptions of theorem~2
are satisfied with $T=\mathbb{R},$ $V=\mathbb{R}^2,$
$E=\mathbb{R},$ $Z=\mathbb{R}^2,$ $\eps_0=0.$ Properties (i) and
(ii) are straightforward and we, therefore, start with (iii).
\begin{itemize}
\item[(iii)] Since $x'_t(t,v,0)\equiv 0$ then
$F'_t(t_0,v_0,0)=-[v_0]_1\sin t_0+[v_0]_2\cos t_0$ which, as it
can be easily verified, equals to 0 if and only if
$F(t_0,v_0,0)\not=0.$
\item[(iv)] The conclusion follows observing that
$
x(t_1,v,\eps)-x(t_2,v,\eps)-x(t_1,v_0,0)-x(t_2,v_0,0)=\eps(y(t_1,v,\eps)-y(t_2,v,\eps)),$
\item[(v)] Follows from the obvious identity
$x(t_1,v_2,\eps)-x(t_1,v_1,\eps)-x(t_2,v_2,0)+x(t_2,v_1,0)=\eps(y(t_1,v_2,\eps)-y(t_1,v_1,\eps)).$
\item[(vi)] $F(t,v,0)=v_1\cos t+v_2\sin t$ and so is continuously
differentiable in $v$ and $t.$
\item[(vii)] Follows from assumption (F).
\item[(viii)] The property holds true with $\{t\}={\rm sign}(t).$
\end{itemize}

Therefore, theorem~2 applies and the function $t=\theta(v,\eps)$
solving (\ref{ps2}) and satisfying (\ref{INE}) exists. Moreover,
$$
  \left\|\left[F'_t(t_0,v_0,0)\right]^{-1}F'_v(t_0,v_0,0)\right\|=
  \left\|
\dfrac{1}{-[v_0]_1\sin t_0+[v_0]_2\cos t_0}(\cos t_0,\sin
t_0)\right\|=R.
$$
To prove (\ref{NV}) we recover that theorem~2 claims that for any
$v\in B_\delta(v_0),$ $\eps\in(0,\delta)$ the function
$F(\cdot,v,\eps)$ has a unique zero in $B_\delta(t_0)$ which is
$\theta(v,\eps).$ On the other hand since $t_0\in(a,b)$ is the
only zero of $F(\cdot,v_0,\eps_0)$ then $\delta>0$ can be
diminished, if necessary, in such a way that $F(\cdot,v,\eps)$
does not vanish on $[a,b]\backslash\{\theta(v,\eps)\}$ for any
$v\in B_\delta(v_0),$ $\eps\in(0,\delta).$ Fix some $v_1,v_2\in
B_\delta(v_0).$ From conclusion (\ref{INE}) of theorem~2 we have
that $\|\theta(v_1,\eps)-\theta(v,\eps)\|\le(R+\Delta)\|v_1-v_2\|$
for any $v\in[v_1,v_2],$ $\eps\in(0,\delta),$ which implies that
$$
  [a,b]\backslash\{\theta(v,\eps)\}\supset[a,b]\backslash\left[\theta(v_1,\eps)-(R+\Delta)\|v_1-v_2\|,\theta(v_1,\eps)+(R+\Delta)\|v_1-v_2\|\right]
$$
for any $v\in[v_1,v_2],$ $\eps\in(0,\delta).$ This finishes the
proof.

\qed

\begin{rem} The conclusion of proposition~1 can not be achieved
with the classical implicit function theorem since nothing can be
said around the continuity of the derivative of $(t,v,\eps)\mapsto
F(t,v,\eps)$ with respect to $t$ unless $\eps=0.$
\end{rem}

\end{document}